\documentstyle[12pt]{article}
\input epsf.tex

\begin{document}

\title
{
About periodic solutions of a planar system modelling neural activity%
\thanks{The work is supported by the Russian Foundation for Fundamental
Research (Grant No. 99-01-00574)}
}

\author{S.~Treskov, E.~Volokitin}
\date{}
\maketitle

\begin{abstract}
We derived explicit symbolic expressions for the first, second, and third
Lyapunov coefficients of the complex focus
of a planar system modelling activity of a neural network.
The analysis of these expressions allowed us to obtain new results about
the number and location of limit cycles in the model.
\end{abstract}

\noindent {\bf Introduction}.
In this paper we study the mathematical model which describe processes
in a neural network consisting from two neurons.
The model was suggested in [1] and is defined as a planar
differential system
\refstepcounter{equation}
%\label{sys1}
$$ \begin{array} {ll}
 \dot u_1=-u_1+q_{11} \varphi(u_1) - q_{12} u_2 + e_1,\\
 \dot u_2=-u_2+ q_{21} \varphi(u_1) +e_2,
\end{array}
\eqno{(1)}
$$
where  $q_{11}, q_{12}, q_{21} >0, e_1, e_2 \in {\bf R}$ are parameters, and
$$
\varphi(u_1)=\frac{1}{1+e^{-4u_1}}.
$$

In [1] it was fulfilled a study of model (1) with
the help of the bifurcation theory and it was listed possible types
of dynamical behavior of the neural network;
with analytic and numerical techniques it was constructed the bifurcation
diagram of system (1) for some parameter values.

In [2] it was proved by analytic methods the correctness of the
bifurcation diagram for some specific cases.

In [2] it was demonstrated that the change of variables
$$
u=u_1, v=\frac{u_2-e_2}{q_{21}},
$$
forms system (1) to the system
\refstepcounter{equation}
%\label{sys2}
$$ \begin{array} {ll}
 \dot u=-u+a \varphi(u) - b v + c,\\
 \dot v=-v + \varphi(u),
\end{array}
\eqno{(2)}
$$
which involves three parameters only
$$
a=q_{11}, b=q_{12}q_{21}, c=e_1-q_{12}e_2, a,b>0, c \in {\bf R}.
$$

We intend to study system (2).

Using our results from [3], we derived analytic formulas for the
first, second, and third Lyapunov coefficients of the complex focus
of system (2). A sophisticated treatment of the formulas
allows us to assert that for some parameter values system (2)
has three concentric limit cycles around a unique steady state,
two of these cycles being stable. Such a phase portrait is absent
in the list of possible phase portraits of (2) from [1,2]
although on our mind it is interesting from a physical point of view.

\quad

\noindent {\bf 1}.
Let us to note the following property of system (2).
If $u(t), v(t)$ is the solution of the system then
$u_1(t)=-u(t), v_1(t)=1-v(t)$ is the solution of the system
$$ \begin{array} {ll}
 \dot u=-u+a \varphi(u) - b v - c + b -a,\\
 \dot v=-v + \varphi(u),
\end{array}
$$
This fact means that if we have a phase portrait of (2)
for parameter values $a=a_1, b=b_1, c=c_1$ then there exists
a symmetric with respect to the point $u=0, v=1/2$ phase portrait of the
system for parameter values $a=a_1, b=b_1, c=c_2=-a_1+b_1-c_1$.
In the parameter space points $(a_1,b_1, c_1)$ and $(a_1, b_1, c_2)$
lie on different sides from the hyperplane $a-b+2á=0$ (on equal distances).
So, it is sufficiently to study system (2) for parameter values
from one of a half-space, for example, from the half-space $a-b+2c \geq 0.$

The mentioned property is obvious from the
bifurcation diagrams constructed with numerical methods in [2].

The authors of [1,2] solved all problems of the number, location,
and types of steady states of system  (2). The system
may have between one and three steady states in the finite part of a
phase plane. If the system has three steady states then one of them is always
a saddle. Two remaining steady states are nodes or foci (stable or
unstable). We distinguish between them a left and right steady state.

The authors of [4] introduced the following method to describe
the phase portraits of system (2).
The symbol $s$ ($u$) denotes a stable (unstable) equilibrium;
the symbol $S$ ($U$), a stable (unstable) limit cycle.
The subscript $1$ refers to the symbol corresponding to the left
equilibrium or to a limit cycle around it;
the subscript $2$, to the symbol corresponding to the right
equilibrium or to a limit cycle around it.
A symbol has no subscripts when the steady state is unique.
The symbols $S$ and $U$ without subscripts refer also to limit cycles
around all three steady states.
If there are several cycles with the same subscript or without subscript,
then they are listed in order from inside.

The authors of [1,2] detected following phase portraits of
system (2).%
\footnote
{The subscript 1 for the number of a phase portrait implies that
a corresponding portrait is symmetric to the portrait with the same
number without the subscript.
}

\refstepcounter{equation}
%\label{sp}
\begin{center}
\begin{tabbing}
1$_{\rm 1}$. \= $s_1U_1u_2S$ \quad \ \ \ \ \=  10$_{\rm 1}$.\= $s_1U_1u_2S_2US$ \quad \=  20$_{\rm 1}$. \= $s_1U_1u_2S_2US$ \= \hfill \kill
1.\> $s_1s_2$ \>  5$_{\rm 1}$.\> $s_1u_2$  \>  9$_{\rm 1}$.\> $s_1U_1u_2S$ \>  \\
2.\> $s_1s_2US$ \>  6.\> $u_1s_2US$  \>  10.\> $u_1s_2S$ \> \\
3.\> $s_1U_1s_2US$ \>  6$_{\rm 1}$.\> $s_1u_2US$  \>  10$_{\rm 1}$.\> $s_1u_2S$ \>  \\
3$_{\rm 1}$.\> $s_1s_2U_2US$ \>  7.\> $s_1U_1s_2S$  \>  11.\> $u_1u_2S$ \> \` (3) \\
4.\> $s_1U_1s_2$ \>  7$_{\rm 1}$.\> $s_1s_2U_2S$  \>  12.\> $s$  \> \\
4$_{\rm 1}$.\> $s_1s_2U_2$ \>  8.\> $s_1U_1s_2U_2S$  \>  13.\> $sUS$  \> \\
5.\> $u_1s_2$ \>  9.\> $u_1s_2U_2S$  \>  14.\> $uS$  \> \\
\end{tabbing}
\end{center}

As a result of our study, we detect one more phase portrait of
system (2) which is described by the sequence
\refstepcounter{equation}
%\label{sp1}
$$
{\rm 15}.\ uSUS
\eqno{(4)}
$$
and corresponds to the situation when in the system two
stable limit cycles coexist which are separated with an unstable cycle.

The detection is due to a more comprehensive analysis of the Lyapunov
coefficients of the complex focus of system (2).
The next section is devoted to this analysis.

\quad

\noindent {\bf 2}.
In [3] we derived symbolic expressions for the first, second, and
third Lyapunov coefficients of the complex focus of the system
$$
\begin{array}{l}
\dot x=y,\\
\dot y=a_{10} x + a_{01} y + a_{20} x^2 + a_{11} x y + a_{02} y^2 + \ldots
\end{array}
$$
in terms of the coefficients $a_{ij}$.

The expressions are rather unwieldy and in general case are not amenable
to theoretical study but in some special cases they may be simplified
noticeably. In particular, if we have a Li\'{e}nard system with the complex
focus in the origin
\refstepcounter{equation}
%\label{sysl}
$$
\begin{array} {ll}
 \dot x=y,\\
 \dot y=p(x)+y q(x) \equiv p_1 x + p_2 x^3 + \ldots + y(q_1 x + q_2 x^2 + \ldots),
\end{array}
\eqno{(5)}
$$
then the expressions for the first three Lyapunov coefficients
with accuracy to a positive factor coincide with ones
\refstepcounter{equation}
%\label{lvl}
$$
\begin{array}{ll}
l_1 = & p_2 q_1 - p_1 q_2,\\
l_2 = & 5(p_2 q_3 - p_3 q_2) + 3(p_4 q_1 - p_1 q_4),\\
l_3 = & 14 p_2 (p_2 q_4 - p_4 q_2) + 21 p_1 (p_3 q_4 - p_4 q_3) +
        35 p_1 (p_5 q_2 - p_2 q_5) + \\
    \ & + 15 p_1 (p_1 q_6 - p_6 q_1).
\end{array}
\eqno{(6)}
$$

Let $(u_0,v_0)$ be a steady state of system (2).
Then we have
$$
u_0-(a-b)\varphi (u_0) -c =0, v_0=\varphi (u_0).
$$

Let us make the following change of variables in system (2)
$$
x=u-u_0, \ y=-u+a \varphi (u) - b v + c.
$$

After it the steady state $(u_0, v_0)$ transfers to the origin.

We obtain a Li\'{e}nard system which is equivalent to system (2)
\refstepcounter{equation}
%\label{sys3l}
$$
\begin{array}{ll}
\dot x=y,\\
\dot y=p(x) + y q(x),
\end{array}
\eqno{(7)}
$$
where
\refstepcounter{equation}
%\label{pq}
$$
\begin{array}{ll}
p(x)= c + (a-b)/(1 + \exp(-4u_0 - 4x)) - u_0 - x,\\
q(x)= -2 + (4 a \exp(-4 u_0 - 4 x)/(1 + \exp(-4 u_0 - 4 x))^2.
\end{array}
\eqno{(8)}
$$

The origin will be the complex focus of system (7), (8)
if $p(0)=q(0)=0, p'(0)<0,$ that is under the fulfillment of the conditions
$$
\begin{array}{ll}
c + (a-b)/(1 + \exp(-4u_0)) - u_0=0,\\
-2 + (4 a \exp(-4 u_0)/(1 + \exp(-4 u_0))^2=0,\\
-1 + 4(a-b) \exp(4 u_0))/(1 + \exp(4 u_0))^2 <0,
\end{array}
$$
which may be written as
\refstepcounter{equation}
%\label{cond}
$$
\begin{array}{ll}
c =(b -a)/(1 + \exp(-4u_0)) + u_0,\\
a = (1 + \exp(4 u_0))^2/(2 \exp(4 u_0)),\\
-1 + 4(a-b) \exp(4 u_0))/(1 + \exp(4 u_0))^2 <0.
\end{array}
\eqno{(9)}
$$

Under conditions (9) system (7), (8)
takes form (5).

Evaluating needed coefficients $p_i, q_i$ and substituting them to (6),
we obtain the expressions for the Lyapunov coefficients of system
(7), (8).

After some transformations we obtain that the Lyapunov coefficients
with accuracy to a positive factor are equal to the expressions

\refstepcounter{equation}
%\label{lv12}
$$
\begin{array}{ll}
l_1 =& 1 - 2 \exp(4 u_0) - (6-8a+8b) \exp(8 u_0) - 2 \exp(12 u_0) +\\
   \ & +\exp(16 u_0),\\
l_2 =& 3 - 72 \exp(4 u_0) + (45 + 8(a-b))\exp(8 u_0) +\\
   \ & (240  - 368(a -b)) \exp(12 u_0) + (45 + 8(a-b))\exp(16 u_0) +\\
   \ & - 72 \exp(20 u_0) + 3 \exp(24 u_0),\\
l_3=&1 -4(29 +a -b)\exp(4u_0) + (717+160(a-b))\exp(8u_0)+\\
   &+16(102 -205(a-b) -6(a-b)^2 )\exp(12u_0) - \\
   &-2(903- 1904(a-b) -384(a-b)^2)\exp(16u_0)+\\
   &+8(651 + 1813(a-b) -1608( a-b)^2 )\exp(20u_0)-\\
   &-2(903 -1904(a-b) -364(a-b)^2)\exp(24u_0)+\\
   &+16(102 -205(a-b) -6(a-b)^2 )\exp(28u_0) + \\
   &+(717+160(a-b))\exp(32u_0) -4(29 +a -b)\exp(36u_0) +\\
   &+\exp(40u_0).

\end{array}
\eqno{(10)}
$$

Conditions (9) define a two-dimensoinal
bifurcation manifold  in the parameter space. This manifold corresponds
to the codim 1 Andronov-Hopf bifurcation. The condition $l_1=0$
from (10) defines in this manifolds a one-dimensional bifurcation
codim 2 manifold corresponding to the degenerate Andronov-Hopf bifurcation
(the Bautin bifurcation, the Takens bifurcation, the vanishing of the first
Lyapunov coefficient).

The last manifold is piecewise smooth and it may be defined as follows

\refstepcounter{equation}
%\label{l1par}
$$
\begin{array}{ll}
a= (1 + \vartheta)^2/(2 \vartheta),\\
b=(1 + \vartheta)^2 (1 + \vartheta^2)/(8 \vartheta^2),\\
c=(1 - 3 \vartheta - 3 \vartheta^2 + \vartheta^3 + 2 \vartheta \ln \vartheta)/(8 \vartheta),\\
0<\vartheta<+\infty, \vartheta \neq 1.
\end{array}
\eqno{(11)}
$$

Coordinates $(u_0,v_0)$ of a  corresponding unique equilibrium may be
obtained from equations $\exp (4 u_0)=\vartheta, v_0=\varphi(u_0).$

If parameters $a, b, c$ are evaluated accordingly to (11)
then system (7), (8) (and then system (2) also)
has the complex focus of the multiplicity 2.
Stability of the focus is defined by the sign of the expression $\bar l_2$
which is the second Lyapunov coefficient $l_2$ from (10) evaluating
accordingly to the (11) and to the equality $\vartheta=\exp(4 u_0).$
$$
\bar l_2=2(1+\vartheta)^2 (1-14 \vartheta + 6 \vartheta^2 -14 \vartheta^3 + \vartheta^4).
$$

The equation $\bar l_2=0$ has two positive roots $\vartheta_1, \vartheta_2$
$$
\begin{array}{ll}
\vartheta_1=(7+3\sqrt{5}+\sqrt{6}\sqrt{15+\sqrt{5}})/2 \approx 13.6349,\\
\vartheta_2=(7+3\sqrt{5}-\sqrt{6}\sqrt{15+\sqrt{5}})/2 \approx .0733414,\\
\end{array}
$$
to which correspond the following parameter values
$a_1 \approx 7.8541, b_1\approx 26.9164,$ $c_1\approx 18.4129$ and
$a_2=a_1, b_2=b_1, c_2 \approx 0.64937,$
$a_1-b_1+2á_1>0, a_2-b_2+2á_2<0.$

The third Lyapunov coefficient is negativ for $(a_1, b_1, c_1)$
and $(a_2, b_2, c_2)$.

For such parameters system (2) has a unique steady state
which is the complex focus of multiplicity 3.

It is well known that in the parameter space to the point corresponding
to the vanishing of the second Lyapunov coefficient is adjacent a region
of parameter values for which the system has three concentric limit
cycles around the equilibrium.

Thus, we proved that there are parameter values $a,b,c,$ for which
system (2) has three limit cycles surrouding a unique steady
state (an unstable focus). Two from these cycles are stable and they
are separated by an unstable cycle.

\quad

\noindent {\bf 3}.
Fig.1 gives a part of a complete bifurcation diagram of system (2)
for a fixed value of parameter $a>a_1.$ The rest of the diagram may be
easily restored from the about mentioned reasons of the symmetry.
The diagram is fairly typical. In particular, it has many common features
with the diagram given in [1,2] for $a<a_1.$ We use the notations
introduced in [1,2] and restrict ourselves to only brief remarks.
More detailed description of the mentioned bifurcations may be found
in [6,7], for example.

The curve $sn_1$ corresponds to the availability in the system of a
a double steady state which by an appropriate parameter perturbation
splits into a saddle and a node which lies lefter then a saddle.
On the curve $sn_1$ there is a point $tb_1,$ which corresponds
to the Bogdanov-Takens bifurcation and divides the curve into two
segments with stable and unstable saddle-nodes.

The line $h_1\ (h_2)$ corresponds to the Andronov-Hopf bifurcation
of the left (right) steady state if parameters are inside the region
of multiplicity of steady states or of a unique steady state.
As we noted above these lines are defined by formulas (9).
The point $dh_1$ on the line $h_1$ corresponds to the degenerate
Andronov-Hopf bifurcation and divides the line into two segments
corresponding to a soft or hard loss in stability. The second
Lyapunov coefficient is positive at the point $dh_1$. Let us note
that this coefficient is negative at the point $dh_1$ on the diagram
presented in [1,2] for $a<a_1$.

The curves $sl_1$ and $sl_2$ correspond to separatrix loops around
left and right steady states. The curve $sl_4$ corresponds to a separatrix
loop around both steady states.

The curve $snpo$ corresponds to a double cycle which is stable from
outside.
The curve $snpo_1$ corresponds to a double cycle which is stable from
inside. The curves $snpo$ and $snpo_1$ coalesce and form a cusp
at the point corresponding to a triple cycle of system (2).

The listed curves divides the parameter space into regions within which
the system has distinct rough phase portraits describing
in (3), (4).

As usual, we give only the scaleless relative positions of the
bifurcation sets. The drawing of these sets to real scale involves technical
dificulties.

The local bifurcation sets $sn_i, h_i, tb_i, dh_i$
can be described by explicit formulas
(see (9), (11), for example).
Therefore, construction of these sets involves no difficulties.
As to sets of nonlocal bifurcations $sl_i, snpo, snpo_1$,
in constructing them we used the results of the bifurcation theory,
and the results of the qualitative theory of ordinary differential
equations, and the results of simulations.

A considerable difference of our bifurcation diagram from the diagram
in [1,2] (for $a<a_1$) is the presence of region 15 and curve
$snpo_1$ whose existence we proved above.

In Fig.2 is shown the phase portrait of system (2)
for parameter values $a=16., b=130., c=111.165$, corresponding
to region 15.%
\footnote
{To confirm the existence of three limit cycles in system (2)
for the specified parameter values we studied the proper Poincar\'{e}
map numerically and established that it has three fixed points.
}

\begin{center}
\quad \epsfbox{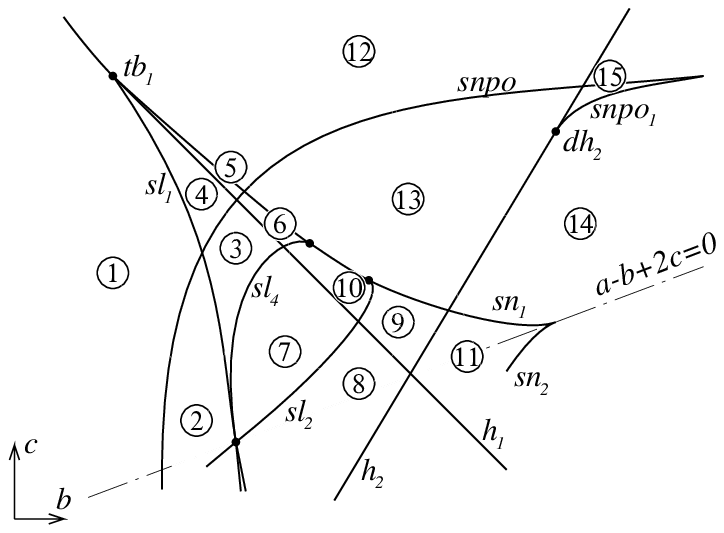} \quad \\
Fig.1\\
\quad \\
\quad \epsfbox{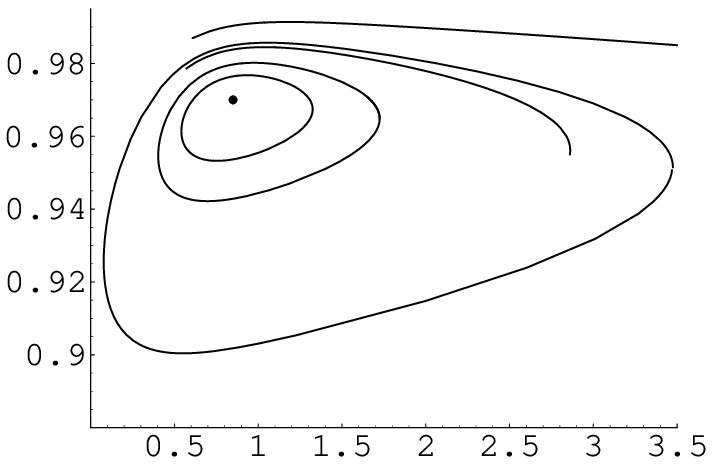} \quad \\
Fig.2
\end{center}

Hence our analysis of formulas for the Lyapunov coefficients for the
system describing the dynamics of a neural network consisting from
two sells enables us to supplement results of other authors.
In particular, to earlier descrided phase portraits we add one more
portrait with three limit cycles around an unstable steady state.
Studing a complete bifurcation diagram, we specify the position
of this portrait in the parameter space.

The results of [1,2] and our results give grounds to advocate
that all possible phase portraits of system (2) are included
in the join of lists (3), (4).

\noindent {\bf 4}.
As we noted above, the bifurcation diagram from [1,2]
is fairly typical and occurs during a study of many mathematical
models, see [8,9] for references.
In particular, a simular diagram was consnructed by us in [10]
where we investigated a planar system modelling a nerve conduction
in the squid giant axon.

{\bf Acknowledgmens.} We thank Professor Fotios Giannakopoulos for materials
which he placed at our disposal.

\end{document}